\documentstyle[11pt,amssymbols]{article}
\font\sym=msbm10 scaled \magstep1
\newcommand{\IP}{{\Bbb P}}
\newcommand{\IR}{{\Bbb R}}
\newcommand{\IH}{{\Bbb H}}
\newcommand{\IC}{\hbox{\sym \char '103}}
\newcommand{\IZ}{\hbox{\sym \char '132}}

\newcommand{\kc}{{\cal C}}
\newcommand{\kd}{{\cal D}}
\newcommand{\kh}{{\cal H}}
\newcommand{\kl}{{\cal L}}
\newcommand{\km}{{\cal M}}
\newcommand{\kn}{{\cal N}}

\newcommand{\kq}{{\cal Q}}
\newcommand{\ko}{{\cal O}}
\newcommand{\kt}{{\cal T}}
\newcommand{\ku}{{\cal U}}
\newcommand{\kv}{{\cal V}}
\newcommand{\ky}{{\cal Y}}
\newcommand{\kx}{{\cal X}}
\newcommand{\kz}{{\cal Z}}

\newcommand{\Pic}{{\rm Pic}}

\newcommand{\Hom}{{\rm Hom }}
\newcommand{\Ext}{{\rm Ext }}
\newcommand{\codim}{{\rm codim}}
\newtheorem{theorem}{Theorem}[section]
\newtheorem{proposition}[theorem]{Proposition}

\newtheorem{lemma}[theorem]{Lemma}
\newtheorem{corollary}[theorem]{Corollary}

\newtheorem{definition}[theorem]{Definition}

\newcommand{\qed}{\hspace*{\fill}\hbox{$\Box$}}
\newcommand{\lra}{\longrightarrow}
\newcommand{\beqn}{\begin{eqnarray}}
\newcommand{\eeqn}{\end{eqnarray}}
\newcommand{\ses}[3]{0\longrightarrow#1\longrightarrow#2\longrightarrow#3
   \longrightarrow0}
\newcommand{\sesq}[5]{0\lra#1\stackrel{#2}{\lra}#3\stackrel{#4}{\lra}#5\lra0}
\newcommand{\rpfeil}[2]{\stackrel{\raisebox{-1mm}{$\scriptstyle#2$}}{{\hbox
to #1mm{\rightarrowfill}}}}
\newcommand{\lpfeil}[2]{\stackrel{\raisebox{-1mm}{$\scriptstyle#2$}}{{\hbox
to #1mm{\leftarrowfill}}}}
\newcommand{\epimorph}{\longrightarrow\!\!\!\!\!\!\!\!\longrightarrow}

\begin{document}
{\parindent0mm{\Large\bf Birational symplectic manifolds and
their deformations}}\\
\bigskip
\bigskip

{\large\bf Daniel Huybrechts}\\
{\small Max-Planck-Institut f\"ur Mathematik, Gottfried-Claren-Str. 26,
53225 Bonn, Germany}
\bigskip

\section{Introduction}
Compact complex manifolds $X^{2n}$ with holonomy group $Sp(n)$ can
algebraically be characterized as simply connected compact K\"ahler
manifolds with a
unique (up to scalars) holomorphic symplectic two-form (\cite{B}).
These manifolds,
which are higher-dimensional analogues of K3 surfaces, are called irreducible
symplectic.

Beauville was able to generalize the local Torelli theorem, one of the
fundamental
results in the theory of K3 surfaces, to all irreducible symplectic
manifolds. His
results show that there exists a (coarse) moduli space $\km$ of
marked
irreducible symplectic manifolds and that the period map
$$P:\km\to\IP(\Gamma\otimes\IC)$$
is \'etale over $Q\subset\IP(\Gamma\otimes\IC)$ -- an open subset of a
quadric defined
by $q(x)=0$ and $q(x+\bar x)>0$.
By definition, a marking is an isomorphism of lattices
$\sigma:H^2(X,\IZ)\cong\Gamma$,
where $H^2(X,\IZ)$ is endowed with the quadratic form defined in \cite{B}
and
$\Gamma$ is a fixed lattice.

For K3 surfaces the moduli space $\km$ consists of two connected
components which
can be identified by $(X,\sigma)\mapsto(X,-\sigma)$. The global Torelli
theorem
for K3 surfaces asserts that the period map $P$ restricted to either of
the two
components, say $\km_0$, is surjective and `almost injective'. More
precisely,
if $(X,\sigma)$ and $(X',\sigma')$ are two points in $P_0^{-1}(x)$, then
$(X,\sigma)$, $(X',\sigma')\in\km_0$ are non-separated and the underlying
$X$ and $X'$
are isomorphic K3 surfaces containing at least one $(-2)$-curve.
Furthermore, for
$x\in Q$ in the complement of the union of countably many proper closed
subsets
the fibre $P_0^{-1}(x)$ is a single point. In short, the failure of the
injectivity
of the period map $P_0$ is due to the non-separatedness of $\km_0$ and
two non-separated points are given by one K3 surface equipped with two
different
markings related by reflections orthogonal to $(-2)$-curves.

In the higher-dimensional situation, the global Torelli theorem does not
hold, i.e.
an isomorphism of Hodge structures $H^2(X,\IZ)\cong H^2(X',\IZ)$ compatible
with the quadratic forms does not imply $X\cong X'$.
In fact, for any two birational irreducible symplectic manifolds $X$ and
$X'$ one
finds markings $\sigma$ and $\sigma'$ such that $P(X,\sigma)=P(X',\sigma')$.
Due to
an example of Debarre, birational $X$ and $X'$ need
not be isomorphic in higher dimensions.

Although, only little evidence can be provided, we cannot resist to formulate
the following (cf. \cite{Mu2}):

{\bf Speculation} {\it (Global Torelli theorem) The period map $P_0$ is
almost
injective, i.e. two points $(X,\sigma)$ and $(X'\sigma')$ in the same fibre
of $P_0$
are non-separated in $\km_0$. In particular, $X$ and $X'$ are birational.}

The birationality of $X$ and $X'$ follows from \cite{MM}.\\
As the known counterexamples to the global Torelli theorem use birational
manifolds
$X$ and $X'$, the following conjecture can be regarded as a weaker version
of this speculation:

{\bf Conjecture} {\it Two irreducible symplectic
manifolds $X$ and $X'$ are birational if and only if they correspond
to non-separated points in
the moduli space.}

This paper proves the conjecture in two fairly general cases.

{\bf Theorem \ref{rho2}} {\it If $X$ and $X'$ are projective irreducible
symplectic
manifolds which are birational and isomorphic in codimension two
(cf. \ref{birat}), then the corresponding points
in the moduli space of symplectic manifolds are non-separated.}

Dropping the assumption on the codimension and the projectivity, but
restricting to
Mukai's elementary transformation, the most explicit birational
correspondence, one can prove

{\bf Theorem \ref{Mukaith}} {\it If $X'$ is the elementary transformation
of an irreducible
symplectic manifold $X$ along a smooth $\IP_N$-bundle of codimension $N$,
then
$X$ and $X'$ correspond to non-separated points in the moduli space.}

Both results combined will be used in Sect. 5 to deduce the conjecture
for projective $X$ and $X'$ and birational correspondences which in
codimension two are
given by elementary transformations (cf. \ref{maincodtwo}).

Unfortunately, only few examples of irreducible symplectic manifolds
are knwon. Higher-dimensional examples were first  described by
Beauville and Fujiki.
Starting with a K3 surface $S$, Beauville showed that the Hilbert schemes
$Hilb^n(S)$ of zero-dimensional subschemes are irreducible symplectic.\\
As shown by Mukai \cite{Mu1}, moduli spaces of stable sheaves on a K3
surface also
admit a (holomorphic) symplectic structure. That these spaces are
irreducible
symplectic, provided they are compact, was shown in \cite{GH} for
the rank two case and
in \cite{OG} in general. The idea in both approaches is to deform the
underlying K3 surface $S$ to a special K3 surface $S_0$, such that the
moduli space
of sheaves on $S_0$ is birational to the Hilbert scheme $Hilb^n(S_0)$.
As the moduli space of sheaves on $S_0$ is a deformation of the moduli
space
of sheaves on $S$, this shows that any smooth moduli space
is deformation equivalent to a manifold which is birational
to an irreducible symplectic manifold. This is enough to conclude that
the moduli spaces of higher rank sheaves are irreducible symplectic.\\
Proving this result \cite{GH}, we observed the following phenomenon.
Let $S$ be a K3 surface and
let $H$ and $H'$ be two different generic polarizations. Then the moduli
spaces
$X:=M_H$ and $X':=M_{H'}$ of $H$-stable, respectively $H'$-stable, sheaves,
which in general are not isomorphic, can be realized as the
special fibres of the same family, i.e.
equipped with appropriate markings they correspond to non-separated points
in the moduli space $\km$. This observation motivated the study of the
general question explained above. Moreover, since the birational
correspondence
between $M_H$ and $M_{H'}$ looks quite similar to the one between moduli
space and
Hilbert scheme on the special K3 surface $S_0$, we conjectured that moduli
spaces
of higher rank sheaves are deformation equivalent to Hilbert schemes
$Hilb^n(S)$.

The general results \ref{rho2} and \ref{Mukaith} do not cover this case,
since
the birational correspondence of moduli space and Hilbert scheme is not
an isomorphism
in codimension two. But using the result of Sect. 5 one can at least
prove the rank two
case.

{\bf Theorem \ref{moddefhilb}} {\it If $S$ is a K3 surface,
$\kq\in\Pic(S)$ indivisible,
$2n:=4c_2-c_1^2(\kq)-6\geq10$ and $H$ a generic polarization,
then the moduli space $M_H(\kq,c_2)$ of $H$-stable rank two sheaves
$E$ with $\det(E)\cong\kq$ and $c_2(E)=c_2$ is deformation equivalent
to $Hilb^n(S)$.}\\

The assumption $2n\geq10$ is a technical condition, whereas the assumption on
the determinant and the polarization is needed to guarantee the smoothness of
the moduli space.
We believe that the same result can be proved for the rank $>2$ moduli
spaces,
as well. As there is evidence that our conjecture holds in general
and that
the higher rank case is an immediate consequence of it, we developed the
necessary
modification only in the rank two case.\\
Due to this result it seems that all known examples of irreducible
symplectic manifolds are either deformation
equivalent to some Hilbert scheme $Hilb^n(S)$, where $S$ is a K3 surface,
or to
a generalized Kummer variety $K^n(A)$, where $A$ is a two-dimensional torus.

{\bf Acknowldegements:} I had valuable and stimulating
discussions with many people.
Especially, I wish to thank A. Beauville, F. Bogomolov, P. Deligne,
B. Fantechi, R. Friedman, S. Keel, and E. Viehweg. I also wish to
thank M. Lehn who has read a first version of the paper.
Part of this work was done during the academic year 1994/95
while I was visiting the Institute for Advanced Study
(Princeton). I was financially supported by a grant from the DFG.
Support and hospitality of the Max-Planck-Institut f\"ur Mathematik (Bonn)
are also gratefully acknowledged.

%%%%%%%%%%%%%%%%%%%%%%%%%%%%%%%%%%%%%%%%%%%%%%%%%%%%%%%%%%%%%%%%%%%%%%%%%%%
\section{Preparations}

\refstepcounter{theorem}\label{sympl}{\bf \thetheorem. Symplectic manifolds.}
A complex manifold $X$ is called {\it symplectic} (in this paper!) if
there exists
a holomorphic two-form $\omega\in H^0(X,\Omega^2_X)$
which is non-degenerate at every point. Note that the existence of $\omega$
implies that the canonical bundle $K_X$ is trivial.
If $X$ is compact, then the symplectic structure is unique if and only if
$h^0(X,\Omega^2_X)=1$.
A simply connected compact K\"ahler manifold with a unqiue symplectic
structure
is called {\it irreducible symplectic}. By \cite{B} $X^{2n}$ is irreduible
symplectic
if and only if its holonomy is $Sp(n)$, i.e. it is irreducible
hyperk\"ahler.\\
For a compact irreducible symplectic K\"ahler
manifold Beauville introduced a
quadratic form on $H^2(X,\IC)$ by
$$\alpha\mapsto\frac{n}{2}\int(\omega\bar\omega)^{n-1}\alpha^2+(1-n)
\int\omega^{n-1}\bar\omega^n\alpha\cdot\int\omega^n\bar\omega^{n-1}\alpha$$
where  $\omega\in H^0(X,\Omega^2_X)=H^{2,0}$ is the symplectic
form. Using Hodge decomposition $\alpha=a\omega+\varphi+b\bar\omega$ with
$\varphi\in H^{1,1}(X)$ and assuming $\int(\omega\bar\omega)^{n}=1$
this form
can be written as
$\alpha\mapsto ab+(n/2)\int(\omega\bar\omega)^{n-1}\varphi$.
It turns out that this form is non-degenerate of index $(3,b_2-3)$.
Moreover, a positive multiple
 of it is integral (cf. \cite{B}, \cite{Fuji}). The unique positive
multiple making it to a primitive integral form
is called the canonical form $q$ on $H^2(X,\IC)$.
Using the weight-two Hodge structure endowed with this quadratic form
Beauville's
local Torelli theorem says that
$\kx_t\mapsto [H^{2,0}(\kx_t)]\in \IP(H^2(X,\IC))$
induces a local isomorphism of the Kuranishi space $Def(X)$ with the
quadric in $\IP(H^2(X,\IC))$ defined by
$q(\alpha)=0$.

\refstepcounter{theorem}\label{birat}{\bf \thetheorem. Birational
symplectic manifolds.}
Let $f:X\to X'$ be a birational map between two compact symplectic
manifolds and assume that the symplectic structure on $X$ is unique.
Then the largest open subset $U\subset X$ where $f$ is regular
satisfies $\codim (X\setminus U)\geq2$.
Moreover, one shows
$f|_U$ is an embedding:
Since $\omega_X$ is unique and $\IC=H^0(X,\Omega_X^2)=H^0(U,\Omega_U^2)$,
the pull-back
$f^*\omega_{X'}$ is a non-trivial multiple of $\omega_X$. Thus $f$ is
quasi-finite on $U$.
Since it is generically one-to-one, it is an embedding.
Note that, as a consequence, the symplectic structure on $X'$ is unique, too.
Thus, if $U\subset X$ and $U'\subset X'$ denote the maximal open subsets
where
$f$ and $f^{-1}$, respectively, are regular, then $U\cong U'$ and
$\codim (X\setminus U)$,
$\codim (X'\setminus U')\geq2$. A birational correspondence is by
definition an
{\it isomorphism in codimension two} if and only if $\codim (X\setminus U)$,
$\codim (X'\setminus U')\geq3$. Recall, that a birational map
between two K3 surfaces
can always be extended to an isomorphism.\\
If $X$ is a projective manifold and $U\subset X$ is an open subset
with $\codim (X\setminus U)\geq2$, then the restriction defines an
isomorphism $\Pic(X)\cong\Pic(U)$.
In particular, for two birational projective manifolds
$X$ and $X'$ with unique symplectic structures one has
$\Pic(X)\cong\Pic(U)\cong\Pic(U')\cong\Pic(X')$. The corresponding
line bundles on $X$ and $X'$ will usually be denoted by
$L$ and $L'$, or $M$ and $M'$. In particular, the Picard numbers
$\rho(X)$ and $\rho(X')$ are equal. Using the exponential sequence
one gets the same result for non-projective
$X$ and $X'$.\\
Frequently, we will use the following result due to  Scheja [S].
If $E$ is a locally
free sheaf on $X$ and $U\subset X$ is an open subset, then
the restriction map $H^i(X,E)\to H^i(U,E|_U)$ is injective for
$i\leq\codim (X\setminus U)-1$ and bijective for
$i\leq\codim (X\setminus U)-2$.
In particular, this can be applied to the line bundles $L$ and $L'$. Thus,
$H^0(X,L)=H^0(U,L|_U)=H^0(U',L'|_{U'})=H^0(X',L')$ and if
$\codim(X'\setminus U')\geq3$
we get $H^1(X,L)\subset H^1(X',L')$.\\
If $X$ and $X'$ are birational irreducible symplectic manifolds,
then there exists an isomorphism between their weight-two Hodge
structures compatible
with the canonical forms $q_X$ and $q_{X'}$ (\cite{Mu2}, \cite{OG}).

\refstepcounter{theorem}\label{defo}{\bf\thetheorem. Deformations.}
Any compact K\"ahler manifold $X$ with trivial canonical bundle $K_X$ has
unobstructed deformations, i.e. the base space of the Kuranishi family
$Def(X)$ is smooth. This is originally due to
Bogomolov, Tian and Todorov  (\cite{Bo2}, \cite{Ti}, \cite{To}).
For an algebraic proof see \cite{R} and \cite{Ka}.\\
If $L$ is a line bundle on $X$, such that the cup-product
$c_1(L):H^1(X,\kt_X)\to H^2(X,\ko_X)$ is surjective, then the deformations
of the pair $(X,L)$ are unobstructed as well. This follows from the fact
that the
infinitesimal deformations of $(X,L)$ are parametrized by $H^1(X,\kd(L))$
and
the obstructions are contained in $H^2(X,\kd(L))$. Here
$\kd(L)$ is the sheaf of differential operators of order $\leq1$ on $L$.
The symbol map
induces an exact sequence
$$\ses{\ko_X}{\kd(L)}{\kt_X}$$
whose boundary map $H^1(X,\kt_X)\to H^2(X,\ko_X)$ is the cup-product with
$c_1(L)$. In particular, $H^2(X,\kd(L))\to H^2(X,\kt_X)$ is injective.
Since $X$ is unobstructed, all obstructions of $(X,L)$ vanish.\\
All this can be applied to irreducible symplectic manifolds.
Using $H^1(X,\kt_X)\cong H^1(X,\Omega_X)$ one finds that $Def(X)$ is
smooth of positive
dimension. Any small deformation of $X$ is again K\"ahler (cf. \cite{KS})
and irreducible symplectic.
In fact, any K\"ahler deformation of $X$ is irreducible symplectic \cite{B}.
Under
the isomorphism $H^1(X,\kt_X)\cong H^1(X,\Omega_X)$ the kernel
of $c_1(L):H^1(X,\kt_X)\to H^2(X,\ko_X)=\IC$ is identified with
the kernel of $q(c_1(L),~~):H^1(X,\Omega_X)\to\IC$ (cf. \cite{B}).
In particular, if
$L$ is non-trivial, then $c_1(L):H^1(X,\kt_X)\to H^2(X,\ko_X)$ is
surjective and thus $Def(X,L)$ is a smooth hypersurface of $Def(X)$.
For the tangent space of $Def(X,L)$ we have
$T_0Def(X,L)\cong H^1(X,\kd(L))\cong\ker(H^1(X,\kt_X)
\rpfeil{5}{c_1(L)}H^2(X,\ko_X))\cong
\ker(H^1(X,\kt_X)\cong H^1(X,\Omega_X)\rpfeil{5}{q(c_1(L),~)}\IC)$.
If $c_1(L)$ and $c_1(M)$ are linearly independent, then the
deformation spaces $Def(X,L)$ and $Def(X,M)$ intersect transversely.

\refstepcounter{theorem}\label{ms}{\bf\thetheorem. Moduli spaces.}
Due to
Beauville's local Torelli theorem one can easily construct a moduli
space $\km$
of marked irreducible symplectic manifolds. Here
a marking consists of an isomorphism of $H^2(X,\IZ)$ with
a fixed lattice compatible with the quadratic form $q$.
As for K3 surfaces the space of marked irreducible symplectic K\"ahler
manifolds is
smooth but non-separated.  In contrast to the K3 surface case, the
moduli space $\km$
is in general not fine. This is due to the fact that higher-dimensional
irreducible symplectic manifolds permit automorphisms inducing the
identity on
$H^2(X,\IZ)$ (cf. \cite{B2}).\\
The quotient of $\km$ by the orthogonal group of $(H^2,q)$ is the moduli
space of unmarked manifolds, but this space
is not expected to have any reasonable analytic
structure. The theme of this paper is to prove statements like:
$X$ and $X'$ correspond to non-separated points in the moduli space.
Here, we usually refer to the moduli space of marked manifolds, though
this distinction does not really matter for our purposes. Explicitly,
this means that there are two one-dimensional deformations $\kx\to S$ and
$\kx'\to S$ ($S$ is smooth),
which are isomorphic over $S\setminus\{0\}$ and the special fibres are
$\kx_0\cong X$ and $\kx'_0\cong X'$.

%%%%%%%%%%%%%%%%%%%%%%%%%%%%%%%%%%%%%%%%%%%%%%%%%%%%%%%%%%%%%%

\section{Elementary transformations}\label{Mukai}
An explicit birational correspondence between two symplectic manifolds
was introduced by Mukai \cite{Mu1}.
We briefly want to recall the construction.\\
Let $X$ be a complex manifold of dimension $2n$ which admits a holomorphic
everywhere non-degenerate two-form $\omega\in H^0(X,\Omega^2_X)$.
Furthermore,
let $P\subset X$ be a closed submanifold which itself is a projective
bundle $P=\IP(F)\rpfeil{5}{\phi}Y$. Here, $F$ is a rank-$(N+1)$ vector
bundle
on the manifold $Y$. Using the symplectic structure one can define the
elementary
transformation $X'$ of $X$ along $P$  as follows.\\
Since a projective space $\IP_N$ does not admit any regular
two-form, the restriction of $\omega$ to any fibre of $\phi$ is trivial.
More is true, the relative tangent bundle $\kt_\phi$ of $\phi$ is
orthogonal to
$\kt_P$ with respect to the restriction of $\omega$, i.e.
$\omega|_P:\kt_\phi\times\kt_P\to\ko_P$
vanishes. Indeed,
this follows from the isomorphism $H^0(Y,\Omega_Y^2)\cong H^0(P,\Omega_P^2)$,
i.e. $\omega|_P$ is the pull-back of a two-form on $Y$. Thus the composition
of $\kt_P\subset\kt_X|_P$ with the isomorphism
$\kt_X|_P\cong\Omega_X|_P$ and the projection
$\Omega_X|_P\epimorph\Omega_P\epimorph\Omega_\phi$
vanishes. Hence $\omega$ induces a vector bundle homomorphism
$\kn_{P/X}\cong\kt_X|_P/\kt_P\to\Omega_\phi$.\\
Now let ${\rm codim}P=N$. Then both vector bundles $\kn_{P/X}$ and
$\Omega_\phi$ are  of
rank $N$ and,
since $\omega$ is non-degenerate, the homomorphism $\kn_{P/X}\to\Omega_\phi$
is an isomorphism.\\
Let $\tilde X\to X$ denote the blow-up of $X$ in $P\subset X$ and let
$D\subset\tilde X$
be the exceptional divisor. The projection $D\to P$ is isomorphic to the
projective bundle $\IP(\kn_{P/X})\cong\IP(\Omega_\phi)\to P$.\\
The natural isomorphism of the incidence
variety $\{(x,H)|x\in H\}\subset\IP_N\times\IP_N^*$ as a projective bundle
over $\IP_N$
with the projective bundle $\IP(\Omega_{\IP_N})\to\IP_N$ can be generalized
to the relative situation, i.e.
there is a canonical embedding
$D=\IP(\Omega_\phi)\subset\IP(F)\times_Y\IP(F^*)$
compatible with the projection to $\IP(F)$. The other projection
$D\to \IP(F^*)$
is a projective bundle as well.
If $\ko_{\tilde X}(D)$ restricts to $\ko(-1)$ on every fibre of
$D\to\IP(F^*)$ then there
exists a blow-down $\tilde X\to X'$ to a smooth manifold $X'$
such that $D\subset \tilde
X$ is the exceptional divisor and $D\to X'$ is the projection
$D\to\IP(F^*)\subset X'$
(cf. \cite{FN}).
Adjunction formula shows that $\ko_{\tilde X}(D)$ indeed
satisfies this condition.
\begin{definition} $X'=elm_PX$ is called the elementary transformation
of the symplectic
mani\-fold $X$ along the projective bundle $P$.
\end{definition}$~$
Mukai also shows that an elementary transformation $elm_PX$ of a
symplectic manifold
$X$ is again symplectic.\\
{\bf Example}\refstepcounter{theorem}\label{ex} {\bf\thetheorem} In
the case of a K3 surface $S$, which
is a two-dimensional symplectic manifold, and a $(-2)$-curve
$P=\IP_1\subset S$ one
obviously has $elm_PS\cong S$. The Hilbert scheme $X:=Hilb^n(S)$, which
is irreducible symplectic, then contains the projective space
$\IP_n\cong S^n(P)= Hilb^n(P)$. The elementary transformation of
$Hilb^n(S)$ along
this projective space is in general not isomorphic to $Hilb^n(S)$.
This is due
to an example of Debarre \cite{Deb}.
Though in his example the K3 surface $S$, and hence $X=Hilb^n(S)$,
is only K\"ahler,
it is expected that one can also find examples
$X\not\cong elm_PX$, where $X$ is projective. Also note that there are
examples where an elementary
transformation of $Hilb^n(S)$ is isomorphic to $Hilb^n(S)$ (\cite{B2}).

The following question was raised in \cite{Mu2}.\\
{\bf Question}\refstepcounter{theorem}\label{Qversion1} {\bf\thetheorem}
{\it Are the symplectic manifolds
$X$ and $X'=elm_PX$ deformation equivalent?}

We want to give an affirmative answer to this question in the case of
compact K\"ahler
manifolds.
\begin{theorem}\label{Mukaith} Let $X$ be a compact symplectic
K\"ahler manifold and let $P\subset X$ be a smooth
$\IP_N$-bundle of codimension $N$. Then there exist two smooth proper
families
$\kx\to S$ and $\kx'\to S$ over a smooth and one-dimensional base $S$,
such that
$\kx$ and $\kx'$ are isomorphic as families over $S\setminus \{0\}$
and the
fibres over $0\in S$ satisfy $\kx_0\cong X$
and $\kx'_0\cong X'\cong elm_PX$.
\end{theorem}

Note that the theorem is in fact stronger than what the original
question suggests.
The theorem shows that $X$ and $X'$ correspond to non-separated points
in the moduli space
of symplectic manifolds.
In particular, one has
\begin{corollary} The higher-weight Hodge structures of $X$ and $elm_PX$
are isomorphic.\qed
\end{corollary}

The following lemma is needed for the proof of the theorem.
Consider a deformation $\kx\to S$ of $X$ and assume that $S$ is smooth and
one-dimensional. Let $v\in H^1(X,\kt_X)$ be its Kodaira-Spencer class, i.e.
$\IC\cdot v$ is the image of the Kodaira-Spencer map
$T_0S\to H^1(X,\kt_X)$. Furthermore, denote by $\bar v\in H^1(X,\Omega_X)$
the image of $v$ under the isomorphism $H^1(X,\kt_X)\cong H^1(X,\Omega_X)$
induced
by the symplectic structure.

\begin{lemma}\label{normaltriv} Assume that $\bar v\in H^1(X,\Omega_X)$ is
a K\"ahler class. Then
the normal bundle $\kn_{P/\kx}$ is isomorphic
to $\phi^*F^*\otimes\ko_\phi(-1)$.
\end{lemma}

{\it Proof:}
We certainly can assume that $Y$ is connected and hence
$H^1(P,\kn_{P/X})\cong
H^1(P,\Omega_\phi)\cong H^0(Y,\ko_Y)\cong\IC$.\\
By construction, the isomorphism $\kn_{P/X}\cong \Omega_\phi$
commutes with the
projections $\kt_X\to\kn_{P/X}$, $\Omega_X\to\Omega_\phi$ and the
symplectic structure
$\kt_X\cong\Omega_X$. In particular, the image $\xi$ of $v$ under
$H^1(X,\kt_X)\to H^1(P,\kn_{P/X})$ is non-zero if and only if
$\bar v$ maps to a non-zero class under
$H^1(X,\Omega_X)\to H^1(P,\Omega_\phi)$. Since $\bar v$ is K\"ahler and
thus its
restriction to the fibres of $\phi$ non-trivial, one concludes that $\xi$
is the
extension class of the unique (up to scalars) non-trivial extension
of $\ko_P$ by $\kn_{P/X}\cong\Omega_\phi$. Thus it is isomorphic
to the
relative Euler sequence
$$\sesq{\Omega_\phi}{}{\phi^*F^*\otimes\ko_\phi(-1)}{}{\ko_P}.$$
Therefore, it suffices to show that $\xi$ is also the extension class
of the canonical  sequence
$$\sesq{\kn_{P/X}}{}{\kn_{P/\kx}}{}{\kn_{X/\kx}|_P},$$
where we use $\kn_{X/\kx}\cong\ko_X$. This follows easily from
the definition of the Kodaira-Spencer class
$v$ as the extension class of
$$\sesq{\kt_X}{}{\kt_\kx|_X}{}{\kn_{X/\kx}}.$$
\qed

{\bf Proof of \ref{Mukaith}:} By \ref{defo}
a one-dimensional deformation $\kx\to S$ of
$X$ such that $\bar v$ is K\"ahler always exists.
Denote the blow-up of $\kx$ in $P$ by
$\tilde \kx\to\kx$. By lemma \ref{normaltriv} the exceptional divisor
${\cal D}\to P$
is isomorphic to the projective bundle $\IP(\phi^*F^*)\to P$. Obviously,
$\IP(\phi^*F^*)\cong\IP(F)\times_Y\IP(F^*)$. Now consider the second
projection
${\cal D}\cong\IP(\phi^* F^*)\to\IP(F^*)$. As before one checks that
${\cal O}_\kx({\cal D})$ restricts to ${\cal O}(-1)$ on every fibre of this
projection, i.e. the condition of the Nakano-Fujiki criterion is satisfied.
Thus $\tilde\kx$ can be blown-down to a smooth manifold $\kx'$ such that
the exceptional divisor ${\cal D}$ is contracted to $\IP(F^  *)$.
By the very construction $\kx'\leftarrow\tilde \kx\to\kx$ is
compatible with $X'\leftarrow\tilde X\to X$, i.e. $\kx'\to S$ is a smooth
proper family, isomorphic to $\kx$ over $S\setminus\{0\}$, and its special
fibre
$\kx'_0$ is isomorphic to $X'$.\qed
\bigskip\\
Note that the two families $\kx$ and $\kx'$ are not isomorphic. In
particular,
one gets the well-known
\begin{corollary} If $X$ is a K3 surface with a $(-2)$-curve $P\subset X$,
then there
exist non-isomorphic families $\kx, \kx'\to S$ which are isomorphic over
$S\setminus\{0\}$
and $\kx_0\cong\kx'_0\cong X$.\qed
\end{corollary}

%%%%%%%%%%%%%%%%%%%%%%%%%%%%%%%%%%%%%%%%%%%%%%%%%%%%%%%%%%%%%%%%%%%%%%%%%%%
\section{Non-separated points in the moduli space}\label{general}

In this section we discuss other situations where birational symplectic
manifolds
present non-separated points in their moduli space.\\
Elementary transformations, dealt with previously, define very explicit
birational
correspondences between symplectic manifolds. But birational correspondences
encountered in the examples are usually more complicated.
This section is devoted to general birational
correspondences. The result is analogous to
\ref{Mukaith}, though we restrict to projective manifolds and
birational correspondences which are isomorphisms in codimension two.
Later (cf. Sect. \ref{codtwo}) the result will be generalized to the case
where in codimension two the birational correspondence is given by an
elementary
transformation.

Let us fix the following notations: $X$ and $X'$ denote irreducible
symplectic
manifolds which are isomorphic on the open sets $U\subset X$ and
$U'\subset X'$ (cf. \ref{birat}).
If $v$ is a class in $H^1(X,\kt_X)$, then the symplectic structure
$\kt_X\cong\Omega_X$
induces a class $\bar v\in H^1(X,\Omega_X)$. The following proposition
does not make any assumptions
either on the projectivity of $X$ or on the codimension of $X\setminus U$.
It is not needed for the proof of the main theorem, but shows how
and to what extent the idea of Sect. 3 works in the general context.
\begin{proposition}\label{ampleKS} Let $S$ be smooth and one-dimensional
and let
$\kx\to S$ and $\kx'\to S$ be deformations of $\kx_0=X$ and $\kx'_0=X'$,
respectively.
If $\kx$ and $\kx'$ are $S$-birational and the Kodaira-Spencer class $v$
of $\kx\to S$
induces a class $\bar v\in H^1(X,\Omega_X)$ which is non-trivial on all
rational
curves in $X\setminus U$, then
$\kx|_{S\setminus\{0\}}\cong_S\kx'|_{S\setminus\{0\}}$
(possibly after shrinking $S$ to an open neighbourhood of $0$).
\end{proposition}

\noindent{\bf Remarks }\refstepcounter{theorem}\label{RMtoampleKS}
{\bf\thetheorem}
{\it i)} $\bar v$ non-trivial on a rational curve means that the pull-back
of $\bar v\in H^2(X,\IC)$ evaluated on the fundamental class
of such a curve is non-trivial.\\
{\it ii)} The condition on $v$ is satisfied if $\bar v$ is contained in
the cone spanned
(over $\IR$) by classes which are ample on $X\setminus U$, e.g. if
$\bar v$ is ample.
Note that the rational curves could be singular and reducible.\\
{\it iii)} Whenever $X$ is projective there are deformations with
Kodaira-Spencer class
$v$ such that $\bar v$ is ample. The problem is to construct $\kx'\to S$
simultaneously.
If the codimensions of $X\setminus U$ and $X'\setminus U'$ are at least
three, then the
isomorphisms
$H^1(X,\kt_X)\cong H^1(U,\kt_U)\cong H^1(U',\kt_{U'})\cong
H^1(X',\kt_{X'})$ suggest that
deformations of $X$ can be related to deformations of $X'$ via
the big open subsets
$U$ and $U'$. I don't know how to make this rigorous. In particular,
it is not clear
to me what deformations of $U$ should really mean.\\
{\it iv)} In the proof of \ref{Mukaith} the family $\kx'\to S$ was
constructed explicitly
from $\kx\to S$ as a blow-up followed by a blow-down. For the general
situation this approach seems to fail.\\

{\bf Proof of \ref{ampleKS}:} If the $S$-birational map $\kx\to\kx'$
does not extend
to an isomorphism $\kx_t\cong\kx'_t$ for generic $t$, then there exists a
surface $\kc$ together with a flat morphism $\kc\to S$ such that:\\
{\it i)} $\kc$ is smooth and irreducible.\\
{\it ii)} For generic $t$ the fibre $\kc_t$ is a disjoint union of
smooth rational
curves.\\
{\it iii)} There exists a finite $S$-morphism $\alpha:\kc\to\kx$ that
maps $\kc_0$
to $X\setminus U$.\\
This follows from resolution of singularities:
By shrinking $S$ we can assume that there is a sequence of monoidal
transformations
$\kz_n\to\kz_{n-1}\to...\to\kz_1\to\kx'$ with smooth centers, which
either dominate $S$ or are contained in the fibre over $0\in S$, and
such that
there exists a morphism $\kz_n\to \kx$ which resolves the birational
map
$\kx\to \kx'$.
If $\kx_t\to\kx'_t$ does not extend to an isomorphism for generic $t$,
then at least one monoidal transformation $\kz_i\to\kz_{i-1}$ with
smooth center
$T_i$ dominating $S$ occurs. Let $i$ be
maximal with this property. Next one finds a morphism
$S'\to T_i$ from a smooth, irreducible curve $S'$ such that the composition
$S'\to T_i\to S$ is finite and smooth over $S\setminus\{0\}$.
Then $\kz_i\times_{\kz_{i-1}}S'\to S'$ is a projective bundle.
Since $i$ is maximal, we have
$(\kz_i\times_{\kz_{i-1}}S')\times_S S\setminus\{0\}\subset\kz_n\times_SS'$,
Now pick a $\IP_1$-bundle contained in
$\kz_i\times_{\kz_{i-1}}S'\to S'$ such that its restriction
to $S'\times_SS\setminus\{0\}$ maps generically finite
to $\kx$ under $\kz_n\to\kx$.
The resolution of the closure of it in $\kz_n$ gives the surface $\kc$.

Now we want to show how one can use the existence of $\kc$ to derive a
contradiction.
First, we claim that the composition
$$T_tS\to H^1(\kx_t,\kt_{\kx_t})\cong
H^1(\kx_t,\Omega_{\kx_t})\rpfeil{5}{\alpha_t^*}H^1(\kc_t,\Omega_{\kc_t})$$
vanishes for generic $t$ (Here, the first map is the Kodaira-Spencer map and
the isomorphism is induced by the symplectic structure on $\kx_t$).
This is a generalization of an argument explained in the proof of
\ref{normaltriv}.
One can either use deformation theory to show that the
existence of $\kc\to S$ implies the vanishing of the obstruction to deform
$\kc_{t\not=0}\to\kx_{t\ne0}$, which in turn gives the desired vanishing,
or one makes this
explicit by the following argument:
Note that we can assume $\kc_{t\not=0}\cong\IP_1$. Then, let
$\kn_t$ be the generalized normal sheaf of $\alpha_t$, i.e. the cokernel
of the injection $\kt_{\kc_t}\to\alpha_t^*\kt_{\kx_t}$.
Since for $t\ne0$ we know
$\kt_{\kc_t}\cong\kt_{\IP_1}$ and $\Hom(\kt_{\IP_1},\Omega_{\IP_1})=0$,
the pull-back of the symplectic structure on $\kx_t$
to $\kc_t$ induces for $t\ne0$ a commutative diagram
$$\begin{array}{ccc}
\alpha_t^*\kt_{\kx_t}&\epimorph&\kn_t\\
\downarrow&&\downarrow\\
\alpha_t^*\Omega_{\kx_t}&\to&\Omega_{\kc_t}\\
\end{array}$$
Thus, in order to show that $T_tS\to H^1(\kc_t,\Omega_{\kc_t})$ is
trivial, it
is enough to prove that
$T_tS\to H^1(\kx_t,\kt_{\kx_t})\to H^1(\IP_1,\alpha_t^*\kt_{\kx_t})\to
H^1(\IP_1,\kn_t)$
vanishes. The image of this map is spanned by the extension class of
$$\ses{\kn_t}{\alpha_t^*(\kt_{\kx}|_{\kx_t})/\kt_{\kc_t}}
{\alpha_t^*(\kt_\kx|_{\kx_t}/\kt_{\kx_t})}$$
(cf. proof of \ref{normaltriv}). Since $\kn_{\kc_t/\kc}$ together
with the natural inclusion
$\kn_{\kc_t/\kc}\subset\alpha_t^*(\kt_\kx|_{\kx_t})/\kt_{\kc_t}$
induced by $\kt_\kc\to\alpha^*\kt_\kx$ splits this sequence, we conclude
that $T_tS\to H^1(\IP_1,\kn_t)$ is trivial.
Hence $T_tS\to H^1(\kc_t,\Omega_{\kc_t})$ is trivial as well.\\
The Kodaira-Spencer map $\kt_S\to R^1\pi_*\kt_{\kx/S}$ composed
with the isomorphism $R^1\pi_*\kt_{\kx/S}\cong R^1\pi_*\Omega_{\kx/S}$
provides a global section of $R^1\pi_*\Omega_{\kx/S}\otimes \Omega_S$.
Trivializing $\kt_S$ we can think of it as an element
in $H^0(S,R^1\pi_*\Omega_{\kx/S})$ or, using Hodge decomposition,
as a $C^\infty$-section of $R^2\pi_*\IC_\kx\otimes\ko_S$.
Moreover, making $S$ small enough we have $R^2\pi_*\IC_\kx\cong
H^2(\kx,\IC)$.
Thus $\kt_S\to R^1\pi_*\kt_{\kx/S}$ induces a $C^\infty$-map
$t\mapsto \bar v_t\in H^2(\kx,\IC)$. Restricting it to  $\kc$ we
get
$\bar w_t\in H^2(\kc,\IC)$. The vanishing we just proved implies
$\langle\bar w_t,[\kc_t]\rangle=0$ for $t\ne0$.
Since also $t\mapsto[\kc_t]\in H^2(\kc,\IC)$ is continous, we can
conclude
$\langle\bar w_0,[\kc_0]\rangle=0$. This contradicts the assumption on
$\bar v|_{X\setminus U}$, since $\kc_0$ as a
degeneration of rational curves is still rational, though singular,
reducible or
even non-reduced \cite{Sa}.\qed

If we in addition assume that $X'$ is projective, then birational
deformations of $X$
and $X'$ can be produced using the following proposition.
\begin{proposition}\label{firstlemma} Suppose $L'\in Pic(X')$ is very ample
and the corresponding line bundle $L\in\Pic(X)$ satisfies $H^1(X,L^n)=0$ for
$n>0$. Let $\kx\to S$ be a deformation of $X=\kx_0$ over a smooth and
one-dimensional base $S$ and assume that there exists a line bundle $\kl$ on
$\kx$ such that
$\kl_0:=\kl|_{\kx_0}\cong L$. Then, replacing $S$ by an open neighbourhood
of $0\in S$
if necessary, there
exists a deformation $\kx'\to S$ of $\kx'_0=X'$ which is $S$-birational
to $\kx$.
\end{proposition}

{\it Proof:} First, shrink $S$ to the open subset of points $t\in S$
such
that $H^1(\kx_t,\kl_t)=0$.
Since $H^1(X,L)=0$, this is an open neighbourhood of $t=0$.
By base change theorem (cf. \cite{Ha}, III. 12.11) $h^0(\kx_t,\kl_t)$
is constant
on $S$, since it can only jump at a point $t$ if $H^1(\kx_t,\kl_t)\not=0$.
Hence $\pi_*\kl$ is locally free on $S$ with fibre
$(\pi_*\kl)(t)=H^0(\kx_t,\kl_t)$.\\
By the very ampleness of $L'$ the base locus $Bs(L)$ of $L$ is contained
in $X\setminus U$ and therefore of codimension at least $2$.
The set $\cup_{t\in S}Bs(\kl_t)$ is a closed subset of $\kx$ and
hence $\codim_{\kx_t} Bs(\kl_t)\geq2$ for $t$ in an open neighbourhood
of $t=0$
(semicontinuity of the fibre dimension).
Since $Bs(\kl_t^n)\subset Bs(\kl_t)$ we can assume that
$\codim_{\kx_t} Bs(\kl_t^n)\geq2$ for all $n>0$ and $t\in S$.\\
The rational maps $\phi_{|\kl_t|}:\kx_t - -\to\IP(H^0(\kx_t,\kl_t)^*)$,
defined by the complete linear system $|\kl_t|$, glue to a rational $S$-map
$\phi:\kx - - \to\IP((\pi_*\kl)^*)$.
Then $\phi$ is regular at all points of
$\kx_t\setminus Bs(\kl_t)$ ($t\in S$).
Let $\kz$ be the scheme-theoretic closure of the graph $\Gamma_\phi$ of
$\phi$ in
$\kx\times_S\IP((\pi_*\kl)^*)$, i.e. the closure of $\Gamma_\phi$ with
the reduced induced structure.\\
The projection $\varphi:\kz\to\kx$ is isomorphic over every point of
$\kx_t\setminus Bs(\kl_t)$, $t\in S$.
Note that a fibre $\kz_t$ of $\kz$ over $t\in S$ does not necessarily
coincide with the closure
of the graph of $\phi_{|\kl_t|}$.
However, since $\kx$ has irreducible fibres and hence $\Gamma_\phi$,
the generic
fibre of $\kz\to S$ is irreducible as well. Thus,
shrinking to an open neighbourhood of $t=0$, we can assume
that $\kz_t$ is irreducible for $t\not=0$. In particular,
$\kz_{t\not=0}$ equals the closure of the graph of $\phi_{|\kl_t|}$ in
$\kx_t\times\IP(H^0(\kx_t,\kl_t)^*)$ at least set-theoretically.
Since $\kz$ is integral, i.e. irreducible and reduced, and $S$ is smooth and
one-dimensional, the dominant morphism $\kz\to S$ is flat (\cite{Ha},
III. 9.7.).\\
Now consider the other projection $\psi:\kz\to\IP((\pi_*\kl)^*)$ and
denote its image by $\kx'\subset\IP((\pi_*\kl)^*)$. Strictly speaking,
$\kx'$ is the scheme-theoretic image of $\psi$ and since $\kz$ is reduced,
this
is the image with the reduced induced structure. Since $\kx'$ then is
integral
and $\kx'\to S$ is dominant, $\kx'$ is flat over $S$.\\
Obviously, $X'$ is contained in $\kx_0'$. To conclude that $X'=\kx_0'$ it
is enough to show that
$h^0(X',\ko(n)|_{X'})\geq h^0(\kx_0',\ko(n)|_{\kx_0'})$ for $n\gg0$, where
$\ko(1)$ is the tautological ample line bundle on $\IP(H^0(\kx_0,L)^*)$.
Since $\ko(1)|_{X'}\cong L'$ and $h^0(X',L'^n)=h^0(X,L^n)$,
this is equivalent to $h^0(X,L^n)\geq h^0(\kx_0',\ko(n)|_{\kx_0'})$ for
$n\gg0$. For any $n$ there exists an open neighbourhood $S_n\subset S$ of
$0\in S$, such that $H^1(\kx_t,\kl_t^n)=0$ for $t\in S_n$. This follows
from the vanishing
of $H^1(X,L^n)$ for all $n$. On the intersection
$\cap S_n\subset S$, which is the complement of countably many points,
all the cohomology groups $H^1(\kx_t,\kl_t^n)$ vanish and therefore
$h^0(\kx_t,\kl_t^n)=h^0(X,L^n)$.
Using this and the flatness of $\kx'\to S$, the inequality
$h^0(X,L^n)\geq h^0(\kx_0',\ko(n)|_{\kx_0'})$ is equivalent to
$h^0(\kx_t,\kl_t^n)\geq h^0(\kx_t',\ko(n)|_{\kx_t'})$ for $n\gg0$.
But the latter can be derived using the composition
$$H^0(\kx_t',\ko(n)|_{\kx_t'})\rpfeil{5}{\psi^*} H^0(\kz_t,\psi^*\ko(n))
\rpfeil{5}{i^*}
H^0(\kx_t\setminus Bs(\kl_t),\kl_t^n)\cong H^0(\kx_t,\kl_t^n).$$
Indeed, $\psi^*$ is injective since $\kz_t\to\kx'_t$ is surjective, and
$i^*$ is injective, since
it is induced by the dense open embedding $\kx_t\setminus
Bs(\kl_t)\subset\kz_t$
($t\not=0$). The last isomorphism is a consequence of
$\codim(\kx_t\setminus Bs(\kl_t))\geq2$.
This shows that $\kx'_0=X'$. Shrinking $S$ further we can also
assume that
$\kx'\to S$ is smooth (\cite{Ha},III. Ex. 10.2).\\
It remains to show  the assertion on the
birationality. Let $\kz^*$ and $\kx'^*$ denote the fibre products
$\kz\times_S(S\setminus\{0\})$ and $\kx'\times_S(S\setminus\{0\})$,
respectively.
Stein factorization decomposes $\psi:\kz^*\to\kx'^*$ into a finite
morphism $f:\ky\to\kx'^*$ and a morphism $\kz^*\to \ky$ with connected
fibres.
One first shows that $f:\ky\to\kx'^*$ is in fact an isomorphism.
Since $f_t:\ky_t\to\kx_t'$ is finite, the line bundle
$f_t^*\ko(1)$ is ample. Thus $f_t^*\ko(n)$ is very ample
for $n\gg0$. In order to prove that $f$ is an isomorphism, it is
therefore
enough to show that
$f_t^*:H^0(\kx_t',\ko(n))\to H^0(\ky_t,f_t^*\ko(n))$ is
surjective. We argue as above: Consider
$$H^0(\kx_t',\ko(n)|_{\kx_t'})\rpfeil{5}{f_t^*} H^0(\ky_t,f_t^*\ko(n))
\hookrightarrow H^0(\kz_t,\psi^*\ko(n))\hookrightarrow H^0(\kx_t,\kl_t^n)$$ and
use
$h^0(\kx'_t,\ko(n)|_{\kx'_t})=h^0(\kx_0',\ko(n)|_{\kx_0'})=h^0(X',L'^n)=h^0(X,L^n)=h^0(\kx_t,\kl_t^n)$
for all $t\in\cap S_n$.
Hence $f_t^*$ is bijective for $t$ in the complement of countably many
points
and therefore $\ky\cong\kx'^*$ after shrinking $S$.
Thus $\kz^*\to\kx'^*$ has connected fibres. On the other hand
$\dim\kz_t=\dim\kx_t=
\dim\kx'_t$. Hence $\kz_t\to\kx_t'$ is birational for $t\not=0$.
\qed

Note that the condition $H^1(X,L^n)=0$ is automatically satisfied if
$\codim (X'\setminus U')\geq3$, i.e. if $X'$ and $X$ are isomorphic
in codimension two. Indeed,
$H^1(X,L^n)\subset H^1(U,L^n|_U)=H^1(U',L'^n_{U'})=H^1(X',L'^n)=0$
by Kodaira vanishing and \cite{Sch}. It is at this point that
the assumption on the codimension of $X\setminus U$ enters.\\
Also note that the existence of $\kl$
implies $q(c_1(L),\bar v)=0$, where $\bar v\in H^1(X,\Omega_X)$ is
induced by the
Kodaira-Spencer class $v\in H^1(X,\kt_X)$ of $\kx\to S$ (cf. \ref{defo}).\\
Next, combining \ref{ampleKS} and \ref{firstlemma} we get

\begin{corollary}\label{orthog} Let $X$ and $X'$ be birational projective
irreducible symplectic manifolds isomorphic in codimension two. Assume
there exists a line bundle $L\in \Pic(X)$ and a class
$\bar v\in H^1(X,\Omega_X)$
such that:\\
-- The induced line bundle $L'\in \Pic(X')$ is ample.\\
-- The restriction of $\bar v$ to any rational curve in $X\setminus U$
is non-trivial.\\
-- $q(c_1(L),\bar v)=0$.\\
Then $X$ and $X'$ correspond to non-separated points in the moduli space.
\end{corollary}

{\it Proof:} By taking a high power of $L$ we can assume that $L'$ is very
ample.
Furthermore, $H^1(X,L^n)=H^1(X',L'^n)=0$ for $n>0$. The deformation space
$Def(X,L)$ of the pair $(X,L)$ is a smooth hypersurface of $Def(X)$.
Since $q(c_1(L),\bar v)=0$ and $T_0Def(X,L)\cong
\ker(H^1(X,\kt_X)\cong H^1(X,\Omega_X)\rpfeil{5}{q(c_1(L),~)}\IC)$
(cf. \ref{defo}), the class $v\in H^1(X,\kt_X)$ is
tangent to $Def(X,L)$. Therefore, there exist a deformation
$\kx\to S$ over a smooth and one-dimensional base $S$
with Kodaira-Spencer class $v$ and a line bundle $\kl$ on
$\kx$ such that $\kl_0\cong L$. Then Proposition \ref{firstlemma} shows
that there
exists a deformation $\kx'\to S$ of $X'$  which is $S$-birational to $\kx$
and we conclude by Proposition \ref{ampleKS}.\qed

\noindent{\bf Remarks }\refstepcounter{theorem}\label{RMtoorthog}
{\bf\thetheorem}
{\it i)} If $\IP_n\cong P\subset X$ is of codimension $n$, then $X$
and $X':=elm_PX$
satisfy the assumptions of the corollary provided they are projective.
Indeed, if $L'\in\Pic(X')$ is ample, then either there
exists an element $\bar v\in H^1(X,\Omega_X)$ orthogonal
to $c_1(L)$ or $X$ and $X'$ are isomorphic. The restriction $\pm\bar v|_P$
is either
ample, hence non-trivial on any rational curve in $P$, or zero.
In the latter case, change $\bar v$ and $L$ by a small rational multiple of
an ample divisor $H$ on $X$. Thus we get $\bar v_1:=\bar v+\beta c_1(H)$
and $L_1:=L+\gamma H$.
By adjusting $\beta$ and $\gamma$ we can assume $q(c_1(L_1'),\bar v_1)=0$
and $L_1'$
ample for small $\gamma$.
Obviously, $\bar v_1|_P\ne0$ and therefore $\bar v_1$ and $L_1$ satisfy
the conditions
of the corollary. Thus \ref{Mukaith} for elementary transformations along
a projective space can be seen as a corollary of \ref{orthog} if
$X$ and $X'$ are projective. Does \ref{orthog} work for general
elementary transformations?\\
{\it ii)}  It is sometimes hard to check if $\bar v$ and $L$ satisfying
the conditions
of \ref{orthog} can be found.
I don't know the answer for the examples discussed in Sect. \ref{appl}.

Using \ref{firstlemma} one can in fact prove corollary \ref{orthog}
without the assumptions on $v$. The proof
relies on the fact that a compact Moishezon K\"ahler manifold is
projective.
It can be used to prove the following
\begin{lemma}\label{rho1} If $X$ and $X'$ are birational compact
irreducible symplectic
K\"ahler manifolds with $\rho(X)=\rho(X')=1$ and
$X'$ is projective, then $X\cong X'$.
\end{lemma}
{\it Proof:} $X$ is K\"ahler and Moishezon, hence projective. Thus,
if $L'$ is the ample
generator of $\Pic(X')$, then $\Pic(X)=\IZ\cdot L$ and either $L$
or $L^*$ is ample. Since
$H^0(X,L^n)=H^0(X',L'^n)\ne0$ for $n\gg0$, one concludes that $L$
is ample and hence
$X\cong X'$.\qed

Note that the isomorphism can be chosen such that it extends the birational
map.

Here now is the main theorem of this paper.
\begin{theorem}\label{rho2} Let $X$ and $X'$ be projective irreducible
symplectic manifolds
which are birational and isomorphic in codimension two.
Then $X$ and $X'$ correspond to non-separated points in their moduli space.
\end{theorem}
{\it Proof:} Assume $X$ and $X'$ are not isomorphic. Then $\rho(X)\geq2$.
Let $L'$ be very ample on $X'$ and let $L$ be the associated line bundle
on $X$. Then $Def(X,L)\subset Def(X)$ is a smooth hypersurface of positive
dimension
$h^{1}(X,\Omega)-1$.
Since $\Pic(X)$ is countable and any line bundle $M\in \Pic(X)$ defines a
smooth
hypersurface $Def(X,M)$ intersecting $Def(X,L)$ transversely if
$M^n\not\cong L^m$ ($n\cdot m\ne0$) (\cite{B} and \ref{defo}), there exists
a generic smooth
and one-dimensional $S\subset Def(X,L)$ such that $S\cap Def(X,M)=\{0\}$
for all line bundles $M$ linearly independent of $L$. Let
$(\kx,\kl)\to S$ be the associated deformation of $(\kx_0,\kl_0)\cong(X,L)$.
Then $\rho(\kx_t)=1$ for general $t\in S$, i.e.
$t$ in the complement of countably many points. Now apply Proposition
\ref{firstlemma}.
We get a deformation $\kx'\to S$ of $\kx'_0\cong X'$ which is
$S$-birational to $\kx$.
Moreover, the proof of \ref{firstlemma} shows that there is a line
bundle $\kl'$ on $\kx'$
such that $\kl'_0\cong L'$. For small $t$ the fibre $\kx_t$ is still
K\"ahler and
$\kl'_t$ is still ample on $\kx'_t$.
Thus the lemma applies and shows $\kx_t\cong\kx'_t$ for general $t$
extending the
$S$-birational map $\kx- -\to\kx'$. Since the set of points $t\in S$,
where
$\kx_t- - \to\kx'_t$ cannot be extended to an isomorphism is closed, we can
shrink $S$ such that $\kx- - \to\kx'$ is an isomorphism over
$S\setminus\{0\}$. \qed

We want to emphasize that the condition on the codimension
of $X\setminus U$ is only needed in order to apply \ref{firstlemma}.
If for the deformation $\kx\to S$ considered in the proof the dimension
$h^0(\kx_t,\kl_t^n)$ does not jump in $t=0$, then the argument goes through.
This will be discussed in length in the next section.\\
As an immediate consequence of the theorem we have the
\begin{corollary}\label{hodgenumbers} If $X$ and $X'$ are as in
theorem \ref{rho2}, then they
are diffeomorphic and their weight-$n$ Hodge structures are isomorphic for all
$n$.\qed
\end{corollary}

%%%%%%%%%%%%%%%%%%%%%%%%%%%%%%%%%%%%%%%%%%%%%%%%%%%%%%%%%
\section{The codimension two case}\label{codtwo}

As before, let $X$ and $X'$ be birational projective irreducible
symplectic manifolds. Let $L'\in\Pic(X')$ be an ample line bundle
and denote by $L\in\Pic(X)$ the corresponding line bundle on $X$.
The assumption on the codimension of $X\setminus U$ in theorem \ref{rho2}
was only needed to ensure $H^1(X,L^n)=0$ for $n>0$. If
${\rm cod}(X\setminus U)=2$,
then $H^1(X,L^n)$ is not necessarily zero. Indeed, consider an elementary
transformation of a four-dimensional manifold $X$ along a projective
plane $\IP_2\subset X$. Then a standard calculation shows $H^1(X,L^n)\ne0$
if $n\geq2$. The vanishing $H^1(X,L^n)=0$ was only needed at one point
in the line
of arguments.
Namely, we used it in proposition
\ref{firstlemma} to conclude that $h^0(\kx_t,\kl_t^n)\equiv const$
for a family $\kx\to S$. One might hope that this holds for another
reason.
Indeed, if $X'=elm_PX$ is an elementary transformation in codimension two
and $\kx\to S$ is as in \ref{Mukaith}, then $h^0(\kx_t,\kl_t^n)\equiv const$.
To prove this use the family $\kx'\to S$ constructed explicitly in the
proof of
\ref{Mukaith} and the equality
$h^0(\kx_t,\kl_t^n)=h^0(\kx'_t,{\kl'_t}^n)\equiv const$,
since $H^1(X',L'^n)=0$. For a general birational correspondence
the situation is more complicated, since we need
$h^0(\kx_t,\kl_t)\equiv const $ in the first place in order to construct
$\kx'\to S$
(cf. \ref{firstlemma}).

First, we will show that under the above assumption ($L'$ ample)
the condition $h^0(\kx_t,\kl^n_t)\equiv const$ holds true infinitesimally,
i.e.
for any deformation $\pi:(\kx,\kl)\to S={\rm Spec} (k[\varepsilon])$
of $(X,L)$ the direct
image $\pi_*\kl$ is locally free. This is not quite enough to prove
\ref{rho2}
in complete generality, but makes it highly plausible.\\
Under an additional assumption
(cf. \ref{ass}) one can in fact prove $h^0(\kx_t,\kl^n_t)\equiv const$.
This is the second goal of the section and the result \ref{maincodtwo}
will be applied
in Sect. \ref{appl} to moduli space and Hilbert scheme on a K3 surface.

Consider $(X,L)$ as above and let $s$ be a global section of $L$. Then
there exists
a Kuranishi space $Def(X,L,s)$ of deformations of the triple $(X,L,s)$
together with the forgetful maps $Def(X,L,s)\to Def(X,L)\to Def(X)$.
The induced map between the tangent spaces of $Def(X,L,s)$ and $Def(X,L)$
is surjective for all $s$ if and only if for any deformation
$\pi:(\kx,\kl)\to {\rm Spec}(k[\varepsilon])$ the direct image $\pi_*\kl$
is locally free.
Therefore, in order to prove that $\pi_*\kl$ is locally free we have
to describe the tangent spaces of $Def(X,L,s)$ and $Def(X,L)$ and the
homomorphism
between them. Note, if one could prove that
$Def(X,L,s)$ is smooth, the infinitesimal result
would immediately imply that $h^0(\kx_t,\kl_t)\equiv const$.\\
We already know $T_0Def(X,L)\cong H^1(X,\kd(L))$ (cf.
\ref{defo}).

\begin{proposition}\label{infconst} {\it i)} The Zariski tangent space
of $Def(X,L,s)$ is naturally isomorphic to the first hypercohomology
of the complex
$$\begin{array}{cccc}
\kd(L,s):&\kd(L)&\rpfeil{5}{s}&L\\
&D&\mapsto&D(s)\\
\end{array}$$
{\it ii)} The map between the Zariski tangent spaces
$\IH^1(X,\kd(L,s))$ and $H^1(X,\kd(L))$ is given by the
$E_1$-spectral sequence relating hypercohomology and cohomology.\\
{\it iii)} If $(X,L)$ and $(X',L')$ are as above, then
$\IH^1(X,\kd(L,s))\to H^1(X,\kd(L))$ is surjective.
\end{proposition}
{\it Proof:} {\it i)} and {\it ii)} are well-known (\cite{We}). For
{\it iii)}
we write down the beginning of the spectral sequence:
$$0\to \IC\to H^0(X,L)\to \IH^1(X,\kd(L,s))\to H^1(X,\kd(L))\rpfeil{5}{s}
H^1(X,L).$$
Therefore, $\IH^1(X,\kd(L,s))\to H^1(X,\kd(L))$ is surjective if and only
if $H^1(X,\kd(L))\rpfeil{5}{s} H^1(X,L)$ vanishes.
Hence, we have to show that the pairing
$$\begin{array}{ccccc}
H^1(X,\kd(L))&\otimes& H^0(X,L)&\to &H^1(X,L)\\
(D&,&s)&\mapsto&D(s)\\
\end{array}$$
is trivial.
Consider the injections $H^1(X,\kd(L))\hookrightarrow H^1(U,\kd(L_U))$ and
$H^1(X,\kd(L'))\hookrightarrow H^1(U',\kd(L'_{U'}))$. It suffices to
show that under the natural isomorphism
$ H^1(U,\kd(L_U))\cong H^1(U',\kd(L'_{U'}))$, given by
$L'_{U'}\cong L_{U}$, the two
spaces are identified. Indeed, if so then the commutative diagram
$$\begin{array}{ccccccc}
H^1(X,\kd(L))&\otimes &H^0(X,L)&\to&H^1(X,L)&\hookrightarrow&H^1(U,L_U)\\
\downarrow\cong&&\downarrow\cong&&&&\downarrow\cong\\
H^1(X',\kd(L'))&\otimes
&H^0(X',L')&\to&H^1(X',L')&\hookrightarrow&H^1(U',L'_{U'})\\
\end{array}$$
and the vanishing $H^1(X',L')=0$ prove the assertion.
In order to compare $H^1(X,\kd(L))$ and $H^1(X',\kd(L'))$ as
subspaces of
$H^1(U,\kd(L_U))\cong H^1(U',\kd(L'_{U'}))$ we make use of the exact sequence
$$\ses{\ko_X}{\kd(L)}{\kt_X}.$$ Its cohomology sequence provides the
short exact
sequence
$$\ses{H^1(X,\kd(L))}{H^1(X,\kt_X)}{H^2(X,\ko_X)}.$$
We first show that
the two subspaces $H^1(X,\kt_X)\hookrightarrow H^1(U,\kt_U)$ and
$H^1(X',\kt_{X'})\hookrightarrow H^1(U',\kt_{U'})$ are identified
under $H^1(U,\kt_U)\cong H^1(U',\kt_{U'})$. Using the symplectic structures
this is equivalent to the analogous statement for $\Omega_X$ and
$\Omega_{X'}$.
Let $X\leftarrow Z\rightarrow X'$ be a smooth resolution of the
birational
correspondence $U\cong U'$. Then
$H^{1,1}(X)\oplus\bigoplus_i\IC\cdot D_i\cong H^{1,1}(Z)
\cong H^{1,1}(X')\oplus\bigoplus_i\IC\cdot D_i$, where the $D_i$'s
are the exceptional
divisors. Since the $D_i$'s are trivial on $U\cong U'\subset Z$, the
induced isomorphism
$H^{1,1}(X)\cong H^{1,1}(X')$ is compatible with restriction.\\ To conclude
the proof
we have to show that under the identification of $H^1(X,\kt_X)$ and
$H^1(X',\kt_{X'})$ as subspaces of $H^1(U,\kt_U)$ the
homomorphisms $c_1(L):H^1(X,\kt_X)\to H^2(X,\ko_X)$ and
$c_1(L'):H^1(X',\kt_{X'})\to H^2(X',\ko_{X'})$ have the same kernel.
Since $\ker(c_1(L)\cdot~)$ is identified with $\ker(q(c_1(L),~~))$
under the isomorphism $H^1(X,\kt_X)\cong H^1(X,\Omega_X)$, this follows
immediately from the fact that $H^2(X,\IC)\cong H^2(X',\IC)$ respects
$q_X$ and $q_{X'}$
(cf. \cite{Mu2},\cite{OG}).\qed

The proposition gives evidence that $h^0(\kx_t,\kl_t)\equiv const$ holds
in general.
In fact, I believe that the same technique should show the vanishing
of the higher obstructions to deform $(X,L,s)$,
but I don't know how to prove this.\\
In the examples it seems as if a birational correspondence between
irreducible symplectic
manifolds might be non-isomorphic in codimension two, but that in such
a case
the birational correspondence is in codimension two given by an elementary
transformation. Thus, it is not completely unlikely, that
the following assumption is always satisfied.  For the
birational correspondence between the moduli space of rank two sheaves
and the Hilbert
scheme this is established in Sect. \ref{appl}.

{\bf Assumption} \refstepcounter{theorem}\label{ass}{\bf \thetheorem}
There exist open subsets
$U\subset V\subset X$ and $U'\subset V'\subset X'$ such that
$\codim(X\setminus V),\codim (X'\setminus V')\geq3$
and $V':=elm_{V\setminus U}V$. In particular, we assume that
$P:=V\setminus U$ is a $\IP_2$-bundle
$\IP(F)\to Y$ over  a smooth not necessarily compact manifold $Y$.
If $X$ and $X'$ are
isomorphic in codimension two we set $U=V$ and $U'=V'$.

We are going to prove \ref{rho2} under this additional assumption
(without using \ref{infconst}).

First note that a modification of the proof of \ref{Mukaith}
immediately yields
\begin{corollary} Assume $X$ and $X'$ satisfy \ref{ass}. If
$\kx\to S$
is a deformation as in the proof of
\ref{Mukaith} (i.e. $\bar v$ is non-trivial on the
fibres of $P\to Y$), then there exists a smooth morphism
$\kv'\to S$ such that $\kv'|_{S\setminus{0}}\cong\kx_{S\setminus\{0\}}$
and
$\kv_0\cong V'$. \qed
\end{corollary}
It can be used to prove
\begin{proposition}\label{cod2th} Let $X$ and $X'$ be as before, in
particular
$L'$ ample, and assume that \ref{ass} is satisfied. If $(\kx,\kl)\to S$ is a
deformation over a smooth and one-dimensional base $S$ such that the class
$\bar v\in H^1(X,\Omega_X)$ associated to the Kodaira-Spencer class is
non-trivial on the fibres of $P\to Y$, then
$h^0(\kx_t,\kl_t)\equiv const$ in a neighbourhood of $t=0$.
\end{proposition}
Since replacing $L'$ by another ample line bundle (if necessary) ensures
that
the generic deformation $\kx\to S$ in $Def(X,L)$ has a Kodaira-Spencer class
$v$ such that $\bar v$ is non-trivial on the fibres of $P \to Y$
(cf. \ref{RMtoorthog}),
the proposition immediately shows
\begin{corollary}\label{maincodtwo}
If $X$ and $X'$ are projective irreducible symplectic
manifolds such that $X'$ is an elementary transformation of $X$ in
codimension
two, i.e. \ref{ass} holds, then $X$ and $X'$ present non-separated points
in the moduli space.\qed
\end{corollary}
{\bf Proof of \ref{cod2th}:} Let $s$ be the local parameter of $S$ at
$0\in S$ and let $S_n$ denote the closed subspace
${\rm Spec}(k[s]/s^{n+1})\subset S$.
Furthermore, let $\kx_n:=\kx\times_S S_n$ and $\kl_n:\kl|_{\kx_n}$.
In order to show
that $h^0(\kx_t,\kl_t)\equiv const$, it suffices to prove that for all
$n$ the restriction $H^0(\kx_n,\kl_n)\to H^0(\kx_{n-1},\kl_{n-1})$ is
urjective.
This will be achieved by comparing it with the analogous restriction maps
for the family $\kv'\to S$. For this purpose we introduce the following
notations. $\ku_n$ denotes the space $(U,\ko_{\kx_n}|_U)$ and is considered
as a deformation of $U$ over $S_n$. Analogously, let
$\kv'_n:=\kv'\times_SS_n$
and $\ku'_n:=(U',\ko_{\kv'_n}|_{U'})$, which is isomorphic to $\ku_n$.
The line bundle $\kl$ induces a line bundle $\kl'$
on $\kv'$. Its restrictions to $\kv'_n$ are denoted by $\kl'_n$.
In particular $\kl'_0$ is isomorphic to $L'|_{V'}$.

First, $H^0(\kv'_n,\kl'_n)\to H^0(\kv'_{n-1},\kl'_{n-1})$ is surjective
for all $n$.
Indeed, using the exact sequence
$$\ses{L'|_{V'}}{\kl'_n}{\kl'_{n-1}}$$
this follows from $H^1(V',L'|_{V'})=H^1(X',L')=0$. Next,
$H^0(\ku'_n,\kl'_n|_{\ku'_n})\to H^0(\ku'_{n-1},\kl'_{n-1}|_{\ku'_{n-1}})$
is surjective
and $H^0(\kv'_n,\kl'_n)\to H^0(\ku'_n,\kl'_n|_{\ku'_n})$ is an isomorphism.
This is
proved by induction starting with $H^0(V', L'|_{V'})=H^0(U',L'|_{U'})$ and
the
commutative diagram
$$\begin{array}{cccccccc}
0\to&H^0(V',L'|_{V'})&\to&H^0(\kv'_n,\kl'_n)&\to&
H^0(\kv'_{n-1},\kl'_{n-1})&\to&0\\
&\downarrow\cong&&\downarrow&&\downarrow\cong&&\\
0\to&H^0(U',L'|_{U'})&\to&H^0(\ku'_n,\kl'_n|_{\ku'_n})&\to&
H^0(\ku'_{n-1},\kl'_{n-1}|_{\ku'_{n-1}})&\to&\\
\end{array}$$
The isomorphism $H^0(\ku'_n,\kl'_n|_{\ku'_n})\cong
H^0(\ku_n,\kl_n|_{\ku_n})$ and
a similar induction argument prove $H^0(\kx_n,\kl_n)\cong H^0(\ku_n,\kl_n)$
and $H^0(\kx_n,\kl_n)\epimorph H^0(\kx_{n-1},\kl_{n-1})$. In the
analogous diagram one  in addition has to use $H^1(X,L)\hookrightarrow
H^1(U,L|_U)$.\qed

%%%%%%%%%%%%%%%%%%%%%%%%%%%%%%%%%%%%%%%%%%%%%%%%%%%%%%%%%%%

\section{Application to moduli spaces of bundles on K3 surfaces}\label{appl}
We briefly recall some facts from \cite{GH} that are necessary for our
purposes.\\
Let $S$ be a K3 surface, let $\kq\in\Pic(S)$ be an indivisible line bundle
and let
$c_2\in\IZ$ such that $2n:=4c_2-c_1^2(\kq)-6\geq10$.
Assume that $H$ is a generic polarization, i.e. an ample line bundle
such that a rank two sheaf $E$ with $\det(E)\cong\kq$ and $c_2(E)=c_2$ is
$H$-semi-stable
if and only if it is $H$-stable. Then the moduli space $M_H(\kq,c_2)$ of
$H$-stable rank-two sheaves with determinant $\kq$ and second Chern number
$c_2$ is
smooth and projective. By \cite{Mu1} the moduli space
$M_H(\kq,c_2)$ admits a symplectic structure.\\
Next, one finds a K3 surface $S_0$ such that $\Pic(S_0)\cong\IZ\cdot H_0$,
where
$H_0$ is ample, and $H_0^2/2+3=n$. In \cite{GH} we showed that under all
these
assumptions the moduli space $M_H(\kq,c_2)$ is deformation equivalent to the
moduli space $M_{H_0}(H_0,n)$ of sheaves on $S_0$. In particular,
$M_H(\kq,c_2)$ is irreducible symplectic if and only if $M_{H_0}(H_0,n)$
is irreducible
symplectic. Moreover, both spaces have the same Hodge numbers.\\
In order to prove that
$M_{H_0}(H_0,n)$ is irreducible symplectic we used Serre correspondence to
relate
this space to the Hilbert scheme $Hilb^n(S_0)$. Roughly, the generic sheaf
$[E]\in M_{H_0}(H_0,n)$ admits exactly one global section and the  zero set
of this
section
defines a point in $Hilb^n(S_0)$. To make this more explicit we consider
the
moduli space $N$ of $H_0$-stable
pairs $(E,s\in H^0(S_0,E))$, such that $\det(E)\cong H_0$ and
$c_2(E)=n$. The parameter in the stability condition for such pairs is chosen
very small and constant. As explained in \cite{GH}
the maps $(E,s)\mapsto Z(s)$ and $(E,s)\mapsto E$ define morphisms
$\varphi:N\to Hilb^n(S_0)$ and $\psi:N\to M_{H_0}(H_0,n)$, respectively.
For the fibers we have
$$\varphi^{-1}(Z)\cong\IP(\Ext^1(I_Z\otimes H_0,\ko_{S_0}))$$
and
$$\psi^{-1}(E)\cong\IP(H^0(S_0,E)).$$
Generically, $h^1(S_0,I_Z\otimes H_0)=1$ and $h^0(S_0,E)=1$.
Thus
$$X:=Hilb^n(S_0)\lpfeil{5}{\varphi} N\rpfeil{5}{\psi} M_{H_0}(H_0,n)=:X'$$
defines
a birational correspondence between irreducible symplectic manifolds.\\
Next, we want to show that $X\lpfeil{5}{\varphi} N\rpfeil{5}{\psi}X'$
satisfies
the assumption \ref{ass}.\\ Using the exact sequence
$$\ses{I_Z\otimes H_0}{H_0}{\ko_Z},$$
the vanishing
$H^1(S_0,H_0)=0$ and $H^2_0/2+3=n$, one shows
$h^1(S_0,I_Z\otimes H_0)=1+h^0(S,I_Z\otimes H_0)$.
Therefore, $f:=\varphi\circ\psi^{-1}$ is regular at points $Z$ which are
not contained
in any divisor $D\in |H_0|$.\\
Let $\kd\to |H_0|$ denote the family of divisors
parametrized by the complete linear system $|H_0|$ and let
$Hilb^n(\kd)\to |H_0|$ be the relative Hilbert scheme.
Then $f$ is regular on the complement $U$ of the image of the natural map
$g: Hilb^n(\kd)\to Hilb^n(S_0)=X$. Since $\dim Hilb^n(\kd)=n+h^0(S_0,H_0)-1=
2n-2=\dim Hilb^n(S_0)-2$, the birational correspondence $f$ is not isomorphic
in codimension two.\\
Let $\kc\to B\subset|H_0|$ denote the family of smooth curves.
The relative Hilbert scheme over $B$ is just the relative symmetric product
$S^n(\kc/B)\to B$, which factorizes naturally through
the relative Picard $\Pic^n(\kc/B)\to B$.\\
The fibre of
the factorization $\phi:S^n(\kc/B)\to \Pic^n(\kc/B)$ over a point
$L\in\Pic^n(\kc_t)$
is naturally isomorphic to $\IP(H^0(\kc_t,L))$.
Note that by Riemann-Roch $\chi(\kc_t,L)=n-H^2_0/2=3$. Hence
$h^0(\kc_t,L)\geq3$.
Let $Y\subset \Pic^n(\kc/B)$ denote the open set of line bundle
$L\in\Pic^n(\kc_t)$
such that $h^0(\kc_t,L)=3$ and let $\phi:P:=\phi^{-1}(Y)\to Y$
be the induced
$\IP_2$-bundle.
\begin{proposition} {\it i)} The morphism $g:Hilb^n(\kd)\to X$ restricted
to $P$ is an embedding.\\
{\it ii)} The union $V$ of $U$ and $g(P)$ is open and
$\codim (X\setminus V)\geq3$.\\
{\it iii)} $V\lpfeil{5}{\varphi}\varphi^{-1}(V)\rpfeil{5}{\psi}
\psi(\varphi^{-1}(V))$
is an elementary transformation along the $\IP_2$-bundle $P$.

\end{proposition}

{\it Proof:} There is a number of little things to check.\\
First, by our
assumption $n\geq5$ we have $H_0^2\geq4$. Thus we can apply a result
of Saint-Donat
(cf. \cite{Mo}) to conclude that $H_0$ is very ample. Hence
$B\subset |H_0|$ is dense.
Moreover, $Hilb^n(\kc)$ is dense in $Hilb^n(\kd)$
(cf. \cite{AIK}, Thm. 5).\\
Next, we show that $Y\subset\Pic^n(\kc/B)$ is non-empty and, therefore,
dense
in $\Pic^n(\kc/B)$.
Indeed, if $x_1,...,x_{n-2}$ are generic points in $S_0$, then there
is exactly one smooth curve $C\in |H_0|$ containing them all.
Let $x_{n-1}$ and $x_{n}$
be two more generic points on $C$ and let $Z:=\{x_1,...,x_n\}$. Then
$h^0(S_0,I_Z\otimes H_0)=1$ and hence $h^1(S_0,I_Z\otimes H_0)=2$. Using the
exact sequence
$$\ses{\ko_{S_0}}{I_Z\otimes H_0}{\ko_C(-Z)\otimes H_0}$$
we get
$h^1(C,\ko_C(-Z)\otimes H_0)=h^1(S_0,I_Z\otimes H_0)+h^2(S_0,\ko_{S_0})=3$
and therefore
$h^0(C,\ko_C(Z))=3$.
Thus the line bundle $L:=\ko_C(Z)$ defines a point in $Y$. Note that one
could invoke
a result by Lazarsfeld \cite{La} to prove $Y\ne\emptyset$. His result
also shows that for the generic curve $\kc_t$ the complemet
of $Y\cap\Pic^n(\kc_t)\subset\Pic^n(\kc_t)$ has at least codimension four.\\
Since $P$ is obviously smooth and any $Z\in {\rm Im}(g)$ satisfies
$h^0(S_0,I_Z\otimes H_0)=1$, the morphism $g$ is an embedding on $P$.\\
By definition $V$ is the intersection of the open set
$\{Z|h^0(S_0,I_Z\otimes H_0)\leq1\}$
and the complement of $g(Hilb^n(\kd)\setminus Hilb^n(\kc))$.
Hence $V$ is open.
The assertion on the
codimension follows from $Y\ne\emptyset$ and the irreducibility of
$Hilb^n(\kd)$ (cf. \cite{AIK}).\\
It remains to prove {\it iii)}. Here we make essential use of the
moduli space $N$.\\
Let $N_P$ denote $\varphi^{-1}(P)$. We first show that $\psi:N_P\to X'$
respects
the projective bundle $\phi:P\to Y$, i.e. a fibre of $\psi$ maps to a
fibre of $\phi$.
Indeed, if $[E]\in X'$ and $s_1,s_2\in H^0(S_0,E)$ are two linearly
independent global sections, then we have a diagram
$$\begin{array}{ccccccc}
&&&\ko_{S_0}&&\ko_{S_0}&\\
&&&s_1\downarrow~~~&&\tilde s_1\downarrow~~~&\\
&\ko_{S_0}&\rpfeil{5}{s_2}&E&\rpfeil{5}{}&I_{Z(s_2)}\otimes H_0&\\
&=\downarrow&&\downarrow&&\downarrow&\\
&\ko_{S_0}&\rpfeil{5}{\tilde s_2}&I_{Z(s_1)}\otimes H_0&\rpfeil{5}{}&\kh&\\
\end{array}$$
Thus $\tilde s_1$ and $\tilde s_2$ vanish along the same curve
$C\in |H_0|$ and
$\kh\cong\ko_C(-Z(s_i))\otimes H_0$ for $i=1,2$.
Hence $\ko_C(Z(s_1))\cong \ko_C(Z(s_2))$, i.e.
$\phi(\varphi(E,s_1))=\phi(\varphi(E,s_2))$.\\
This reduces assertion {\it iii)} to the following problem. Let
$L\in\Pic^n(\kc_t)\cap Y$, let $P_L:=\IP(H^0(\kc_t,L))\cong\IP_2$
and let $N_{P_L}:=\varphi^{-1}(P_L)$, which is a $\IP_1$-bundle over $\IP_2$.
Identify
$P_L\lpfeil{5}{}N_{P_L}\rpfeil{5}{}\psi(N_{P_L})$ with
$\IP_2\lpfeil{5}{}\IP(\Omega_{\IP_2})\rpfeil{5}{}\IP_2^*$ !!\\
The argument goes as follows. Any point $Z\in P_L$ gives an exact sequence
$$\begin{array}{ccc}
\ses{\ko_{S_0}}{I_Z\otimes H_0}{&\ko_{\kc_t}(-Z)\otimes H_0 &}\\
&\cong L^*\otimes K_{\kc_t}&\\
\end{array}$$
Now use the canonical isomorphisms
$$\IP(H^0(\kc_t,L))\cong \IP(H^1(\kc_t,L^*\otimes K_{\kc_t})^*)\cong
\IP(\Ext^1(L^*\otimes K_{\kc_t},\ko_{S_0}))$$
to obtain the exact sequence
$$\ses{q^*\ko_{S_0}\otimes p^*\ko_{P_L}(1)}{I_{\kz}\otimes
q^*H_0\otimes p^*\ko_{P_L}(a)}
{q^*(L^*\otimes K_{\kc_t})},$$
where $q$ and $p$ are the two projections from
$S_0\times P_L$ and $I_\kz$ is the ideal sheaf of the universal subscheme
$\kz\subset S_0\times P_L$.
By restricting to $\{x\}\times P_L$, where $x\in S_0\setminus \kc_t$,
we deduce
$a=1$. The push-forward under $p$ induces the exact sequence
$$\ses{R^1p_*(I_\kz\otimes q^*H_0)\otimes\ko_{P_L}(1)}
{H^1(\kc_t,L^*\otimes K_{\kc_t})\otimes \ko_{P_L}}{\ko_{P_L}(1)}.$$
Hence $R^1p_*(I_\kz\otimes q^*H_0)\cong \Omega_{P_L}$.
It is straightforward to identify $N_{P_L}\to P_L$ with
$\IP(R^1p_*(I_\kz\otimes q^*H_0)^*)$. Thus
$(N_{P_L}\to P_L)\cong (\IP(\kt_{\IP_2})\to\IP_2)\cong
(\IP(\Omega_{\IP_2})\to\IP_2)$.\\
It remains to show that $\varphi:P_E:=\IP(H^0(S_0,E))=\psi^{-1}(E)\to P_L$
is a linear embedding.
On $P_E$ we have
$$\ses{\ko}{q^*E\otimes p^*\ko_{P_E}(1)}{(1\times\varphi)^*
(I_\kz\otimes q^*H_0)\otimes p^*\ko_{P_E}(a)},$$
where by abuse of notation $q$ and $p$ are again the projections from
$S_0\times P_E$.
As above one finds $a=2$. Taking direct images we obtain
$$\ses{H^1(S_0,E)\otimes\ko_{P_E}(1)}{\varphi^*(R^1p_*(I_\kz\otimes
q^*H_0))\otimes \ko_{P_E}(2)}{\ko_{P_E}},$$
i.e. $\ses{\ko_{P_E}(1)}{\varphi^*\Omega_{P_L}\otimes\ko_{P_E}(2)}
{\ko_{P_E}}$. Thus
$\varphi^*\ko_{P_L}(1)\cong\ko_{P_E}(1)$.
\qed

{\bf Remark:} The identification $N_{P_L}\cong\IP(\kt_{P_L})$ seems
to show that
the birational correspondence described by $N$ is not some kind of ``nested
elementary transformation'': It is only in the codimension two case where
one has $\IP(\Omega_{P_L})\cong\IP(\kt_{P_L})$.

Corollary \ref{maincodtwo} now immediately implies
\begin{corollary} If $S_0$ is a K3 surface with $\Pic(S_0)=\IZ\cdot H_0$ and
$H_0^2\geq4$, then $M_{H_0}(H_0,n)$ and $Hilb^n(S_0)$ correspond to
non-separated points
in the moduli space of symplectic manifolds ($n=H_0^2/2+3$).\qed
\end{corollary}
Thus we can conclude
\begin{theorem}\label{moddefhilb} If $S$ is a K3 surface, $\kq\in\Pic(S)$
indivisible,
$2n:=4c_2-c_1^2(\kq)-6\geq10$ and $H$ a generic polarization,
then the moduli space $M_H(\kq,c_2)$ of $H$-stable rank two sheaves
$E$ with $\det(E)\cong\kq$ and $c_2(E)=c_2$ is deformation equivalent
to $Hilb^n(S)$.\qed
\end{theorem}
Note that in particular moduli space and Hilbert scheme are just different
complex structures on the same differentiable manifold.\\
{\bf Remark:} O'Grady works instead of $S_0$ with an elliptic surface and
shows
that every moduli space is deformation equivalent to a moduli space
on an elliptic surface \cite{OG}. The birational correpondence between
moduli space and Hilbert scheme on the elliptic surface is again
given by Serre correspondence. The picture there is slightly more complicated
than what we have encountered above. Nevertheless I believe, that
also in his situation
the assumptions \ref{ass} are satisfied and that moduli space and Hilbert
scheme
are deformation equivalent rank$>2$ as well.

%%%%%%%%%%%%%%%%%%%%%%%%%%%%%%%%%%%%%%%%%%%%%%%%%%%%%%%%%%%%%%

e-mail: huybrech@mpim-bonn.mpg.de
\end{document}